\documentstyle[11pt]{article}
\begin{document}
\setlength{\textheight}{574pt}
\setlength{\textwidth}{432pt}
\setlength{\oddsidemargin}{18pt}
\setlength{\topmargin}{14pt}
\setlength{\evensidemargin}{18pt}
\newtheorem{theorem}{Theorem}[section]
\newtheorem{lemma}{Lemma}[section]
\newtheorem{corollary}{Corollary}[section]
\newtheorem{conjecture}{Conjecture}
\newtheorem{remark}{Remark}[section]
\newtheorem{definition}{Definition}[section]
\newtheorem{problem}{Problem}
\newtheorem{example}{Example}
\newtheorem{proposition}{Proposition}[section]
\title{{\bf DEFORMATION INVARIANCE OF PLURIGENERA}}
\date{February 5, 2001}
\author{Hajime Tsuji}
\maketitle
\begin{abstract}
We prove the invariance of plurigenera under smooth projective 
deformations in full generality.
MSC32J25
\end{abstract}
\section{Introduction}
Let $X$ be a smooth projective variety and let $K_{X}$ 
be the canonical bundle of $X$.
The canonical ring 
\[
R(X,K_{X}):= \oplus_{m\geq 0}H^{0}(X,{\cal O}_{X}(mK_{X}))
\]
is a basic birational invariant of $X$.  
For every positive integer $m$, 
the $m$-th plurigenus $P_{m}(X)$ is defined by 
\[
P_{m}(X) : = \dim H^{0}(X,{\cal O}_{X}(mK_{X})). 
\]
Recently Y.-T. Siu (\cite{si}) proved that for every $m\geq 1$, $P_{m}(X)$ 
is invariant under smooth projective deformations, 
if all the fibers are of general type. 
This result has been slightly generalized by \cite{ka,nak}.

In this paper, we shall prove the invariance of plurigenera
under smooth projective deformations in full generality. 
\begin{theorem}
Let $\pi : X \longrightarrow \Delta$ be a 
smooth projective family of smooth projective varieties
over the unit open disk. 

Then for every positive integer $m$,  
the $m$-th plurigenus $P_{m}(X_{t}) (X_{t} : = \pi^{-1}(t))$
 is independent of $t\in \Delta$. 
\end{theorem}
In the course of the proof of Theorem 1.1, we also prove 
the following theorem. 

\begin{theorem}
Let $\pi : X \longrightarrow \Delta$ be a smooth projective 
family.  Suppose that $K_{X_{0}}$ is pseudoeffective.
Then $K_{X_{t}}$ is pseudoeffective for every $t\in \Delta$. 
\end{theorem}

For the proof of Theoem 1.1, the central problem is the existence of singular hemitian metric $h$ on $K_{X}$ such that the curvature $\Theta_{h}$ is semipositive and 
$h\mid_{X_{t}}$ is an AZD (cf. Definition 2.3) of $K_{X_{t}}$ for every 
$t\in \Delta$. 
As soon as we construct such a metric $h$,  by the $L^{2}$-extension theorem (\cite{o-t}) implies the invariance of 
the plurigenera. 
Here the key point is that the $L^{2}$-extension theorem requires the semipositivity of 
the curvature of the singular hermitian metric,  
but it does not require any strict positivity of the curvature. 
In this sense the $L^{2}$-extension theorem is similar 
to Koll\'{a}r's vanishing theorem (\cite{ko}) in 
algebraic context.
 
The construction of the metric $h$ above consists of the inductive 
estimates of singular hermitian metrics using Bergman kernels. 
This is more straightforward than the inductive comparison of 
multiplier ideal sheaves as in \cite[p. 670, Proposition 5]{si}. 
In this sense the proof is {\bf quantitative} and not qualitative. 
The essential idea of the proof is the {\bf dynamical construction}
of an AZD of the canonical line bundle of 
a smooth projective variety with pseudoeffective canonical line bundle.
This construction works only for pseudoeffective canonical line bundles.
This clarifies  why the canonical line bundle is special.
Again the key ingredient of the proof here is  
the $L^{2}$-extension theorem of holomorphic sections  \cite{o-t,m,o}.
Here  we use the fact that the operator norm of the interpolation 
operator is bounded from above by a {\bf universal constant} 
in an essentialy way. 

Hence the proof here is analytic in nature. 

The author would like to express his hearty thanks to Professor T. Ohsawa 
for his interest in this work and his encouragement. 

\section{Preliminaries}
\subsection{Multiplier ideal sheaves}
In this subsection $L$ will denote a holomorphic line bundle on a complex manifold $M$. 
\begin{definition}
A  singular hermitian metric $h$ on $L$ is given by
\[
h = e^{-\varphi}\cdot h_{0},
\]
where $h_{0}$ is a $C^{\infty}$-hermitian metric on $L$ and 
$\varphi\in L^{1}_{loc}(M)$ is an arbitrary function on $M$.
We call $\varphi$ a  weight function of $h$.
\end{definition}
The curvature current $\Theta_{h}$ of the singular hermitian line
bundle $(L,h)$ is defined by
\[
\Theta_{h} := \Theta_{h_{0}} + \sqrt{-1}\partial\bar{\partial}\varphi ,
\]
where $\partial\bar{\partial}$ is taken in the sense of a current.
The $L^{2}$-sheaf ${\cal L}^{2}(L,h)$ of the singular hermitian
line bundle $(L,h)$ is defined by
\[
{\cal L}^{2}(L,h) := \{ \sigma\in\Gamma (U,{\cal O}_{M}(L))\mid 
\, h(\sigma ,\sigma )\in L^{1}_{loc}(U)\} ,
\]
where $U$ runs over the  open subsets of $M$.
In this case there exists an ideal sheaf ${\cal I}(h)$ such that
\[
{\cal L}^{2}(L,h) = {\cal O}_{M}(L)\otimes {\cal I}(h)
\]
holds.  We call ${\cal I}(h)$ the {\bf multiplier ideal sheaf} of $(L,h)$.
If we write $h$ as 
\[
h = e^{-\varphi}\cdot h_{0},
\]
where $h_{0}$ is a $C^{\infty}$ hermitian metric on $L$ and 
$\varphi\in L^{1}_{loc}(M)$ is the weight function, we see that
\[
{\cal I}(h) = {\cal L}^{2}({\cal O}_{M},e^{-\varphi})
\]
holds.
For $\varphi\in L^{1}_{loc}(M)$ we define the multiplier ideal sheaf of $\varphi$ by 
\[
{\cal I}(\varphi ) := {\cal L}^{2}({\cal O}_{M},e^{-\varphi}).
\]
Similarly we define 
\[
{\cal I}_{\infty}(h) := {\cal L}^{\infty}({\cal O}_{M},e^{-\varphi})
\]
and call it the $L^{\infty}$-{\bf multiplier ideal sheaf} of $(L,h)$.

\begin{definition}
$L$ is said to be pseudoeffective, if there exists 
a singular hermitian metric $h$ on $L$ such that 
the curvature current 
$\Theta_{h}$ is a closed positive current.

Also a singular hermitian line bundle $(L,h)$ is said to be pseudoeffective, 
if the curvature current $\Theta_{h}$ is a closed positive current.
\end{definition}

\subsection{Analytic Zariski decompositions}
In this subsection we shall introduce the notion of analytic Zariski decompositions. 
By using analytic Zariski decompositions, we can handle  big line bundles
like  nef and big line bundles.
\begin{definition}
Let $M$ be a compact complex manifold and let $L$ be a holomorphic line bundle
on $M$.  A singular hermitian metric $h$ on $L$ is said to be 
an analytic Zariski decomposition, if the followings hold.
\begin{enumerate}
\item $\Theta_{h}$ is a closed positive current,
\item for every $m\geq 0$, the natural inclusion
\[
H^{0}(M,{\cal O}_{M}(mL)\otimes{\cal I}(h^{m}))\rightarrow
H^{0}(M,{\cal O}_{M}(mL))
\]
is an isomorphim.
\end{enumerate}
\end{definition}
\begin{remark} If an AZD exists on a line bundle $L$ on a smooth projective
variety $M$, $L$ is pseudoeffective by the condition 1 above.
\end{remark}

\begin{theorem}(\cite{tu,tu2})
 Let $L$ be a big line  bundle on a smooth projective variety
$M$.  Then $L$ has an AZD. 
\end{theorem}
As for the existence for general pseudoeffective line bundles, 
now we have the following theorem.
\begin{theorem}(\cite[Theorem 1.5]{d-p-s})
Let $X$ be a smooth projective variety and let $L$ be a pseudoeffective 
line bundle on $X$.  Then $L$ has an AZD.
\end{theorem}
{\bf Proof of Theorem 2.2}.
Although the proof is in \cite{d-p-s}, 
we shall give a proof here, because we shall use it afterward. 

 Let  $h_{0}$ be a fixed $C^{\infty}$-hermitian metric on $L$.
Let $E$ be the set of singular hermitian metric on $L$ defined by
\[
E = \{ h ; h : \mbox{lowersemicontinuous singular hermitian metric on $L$}, 
\]
\[
\hspace{70mm}\Theta_{h}\,
\mbox{is positive}, \frac{h}{h_{0}}\geq 1 \}.
\]
Since $L$ is pseudoeffective, $E$ is nonempty.
We set 
\[
h_{L} = h_{0}\cdot\inf_{h\in E}\frac{h}{h_{0}},
\]
where the infimum is taken pointwise. 
The supremum of a family of plurisubharmonic functions 
uniformly bounded from above is known to be again plurisubharmonic, 
if we modify the supremum on a set of measure $0$(i.e., if we take the uppersemicontinuous envelope) by the following theorem of P. Lelong.

\begin{theorem}(\cite[p.26, Theorem 5]{l})
Let $\{\varphi_{t}\}_{t\in T}$ be a family of plurisubharmonic functions  
on a domain $\Omega$ 
which is uniformly bounded from above on every compact subset of $\Omega$.
Then $\psi = \sup_{t\in T}\varphi_{t}$ has a minimum 
uppersemicontinuous majorant $\psi^{*}$  which is plurisubharmonic.
We call $\psi^{*}$ the uppersemicontinuous envelope of $\psi$. 
\end{theorem}
\begin{remark} In the above theorem the equality 
$\psi = \psi^{*}$ holds outside of a set of measure $0$(cf.\cite[p.29]{l}). 
\end{remark}

By Theorem 2.3,we see that $h_{L}$ is also a 
singular hermitian metric on $L$ with $\Theta_{h}\geq 0$.
Suppose that there exists a nontrivial section 
$\sigma\in \Gamma (X,{\cal O}_{X}(mL))$ for some $m$ (otherwise the 
second condition in Definition 3.1 is empty).
We note that  
\[
\frac{1}{\mid\sigma\mid^{\frac{2}{m}}} 
\]
gives the weight of a singular hermitian metric on $L$ with curvature 
$2\pi m^{-1}(\sigma )$, where $(\sigma )$ is the current of integration
along the zero set of $\sigma$. 
By the construction we see that there exists a positive constant 
$c$ such that  
\[
\frac{h_{0}}{\mid\sigma\mid^{\frac{2}{m}}} \geq c\cdot h_{L}
\]
holds. 
Hence
\[
\sigma \in H^{0}(X,{\cal O}_{X}(mL)\otimes{\cal I}_{\infty}(h_{L}^{m}))
\]
holds.  
Hence in praticular
\[
\sigma \in H^{0}(X,{\cal O}_{X}(mL)\otimes{\cal I}(h_{L}^{m}))
\]
holds.  
 This means that $h_{L}$ is an AZD of $L$. 
\vspace{10mm} {\bf Q.E.D.} 
\begin{remark}
By the above proof we have that for the AZD $h_{L}$ constructed 
as above
\[
H^{0}(X,{\cal O}_{X}(mL)\otimes{\cal I}_{\infty}(h_{L}^{m}))
\simeq 
H^{0}(X,{\cal O}_{X}(mL))
\]
holds for every $m$. 
\end{remark}
It is easy to see that the multiplier ideal sheaves 
of $h_{L}^{m}(m\geq 1)$ constructed in the proof of
 Theorem 2.2 are independent of 
the choice of the $C^{\infty}$-hermitian metric $h_{0}$.
We call the AZD constructed as in the proof of Theorem 2.2  {\bf a canonical 
AZD} of $L$. 

\subsection{$L^{2}$-extension theorem}
\begin{theorem}(\cite[p.200, Theorem]{o-t})
Let $X$ be a Stein manifold of dimension $n$, $\psi$ a plurisubharmonic 
function on $X$ and $s$ a holomorphic function on $X$ such that $ds\neq 0$ 
on every branch of $s^{-1}(0)$.
We put $Y:= s^{-1}(0)$ and 
$Y_{0} :- \{ x\in Y; ds(x)\neq 0\}$.
Let $g$ be a holomorphic $(n-1)$-form on $Y_{0}$ 
with 
\[
c_{n-1}\int_{Y_{0}}e^{-\psi}g\wedge\bar{g} < \infty ,
\]
where $c_{k}= (-1)^{\frac{k(k-1)}{2}}(\sqrt{-1})^{k}$. 
Then there exists a holomorphic $n$-form $G$ on 
$X$ such that 
\[
G(x) = g(x)\wedge ds(x) 
\]
on $Y_{0}$ and 
\[
c_{n}\int_{X}e^{-\psi}(1+\mid s\mid^{2})^{-2}G\wedge\bar{G} 
\leq 1620\pi c_{n-1}\int_{Y_{0}}e^{-\psi}g\wedge\bar{g}. 
\]
\end{theorem}
\section{Dynamical construction of an AZD}
Let $X$ be a smooth projective variety and let 
$K_{X}$ be the canonical line bundle of $X$. 
Let $n$ denote the dimension of $X$.

Let $A$ be a sufficiently ample line bundle on $X$ 
such that for every pseudoeffective singular hermitian 
line bundle $(L,h_{L})$
\[
{\cal O}_{X}(A+L)\otimes{\cal I}(h_{L})
\]
and 
\[
{\cal O}_{X}(K_{X}+A+L)\otimes{\cal I}(h_{L})\\
\]
are globally generated. 
This is possible by \cite[p. 667, Proposition 1]{si}. 
Let $h_{A}$ be a $C^{\infty}$ hermitian metric on $A$
 with strictly positive curvature. 

{\bf Hereafter we shall assume that} $K_{X}$ {\bf is pseudoeffective}. 
For $m\geq 0$, let $h_{m}$ be the singular hermitian metrics
on $A + mK_{X}$ constructed as follows. 
Let $h_{0}$ be a $C^{\infty}$-hermitian metric $h_{A}$ on 
$A$ with strictly positive curvature. 
Suppose that $h_{m-1} (m\geq 1)$ has been constructed. 
Let $\{ \sigma_{0}^{(m)},\ldots ,\sigma_{N(m)}^{(m)}\}$
be an orthonormal basis of 
$H^{0}(X,{\cal O}_{X}(A+mK_{X})\otimes{\cal I}(h_{m-1}))$ 
with respect to the inner product :
\begin{eqnarray*}
(\sigma ,\sigma^{\prime}) & := &
\int_{X}\sigma\cdot \overline{\sigma^{\prime}}\cdot h_{m-1} \\
 &= & \int_{X}\sigma\cdot \overline{\sigma^{\prime}}\cdot (h_{m-1}\otimes (dV)^{-1})\cdot dV,   
\end{eqnarray*}
where $dV$ is an arbitrary nowhere degenerate $C^{\infty}$ volume form on $X$.
We set 
\[
K_{m} := \sum_{i=0}^{m}\mid \sigma_{i}^{(m)}\mid^{2},
\]
where $\mid \sigma_{i}^{(m)}\mid^{2}$ 
denotes $\sigma_{i}^{(m)}\cdot \overline{\sigma_{i}^{(m)}}$.
We call $K_{m}$ the {\bf Bergman kernel} of $A + mK_{X}$ with respect 
to $h_{m-1}\otimes (dV)^{-1}$ and $dV$.
Clearly it is independent of the choice of the orthonormal basis. 
And we define the singular hermitian metric $h_{m}$ 
on $A + mK_{X}$ by  
\[
h_{m} := K_{m}^{-1}.
\]
It is clear that $K_{m}$ has semipositive curvature 
in the sense of currents. 
We note that for every $x\in X$
\[
K_{m}(x) = \sup \{ \mid\sigma\mid^{2}(x) ;  
\sigma\in \Gamma (X,{\cal O}_{X}(A+mK_{X})), 
\int_{X}h_{m-1}\cdot \mid\sigma\mid^{2} = 1\}
\]
holds by definition (cf. \cite[p.46, Proposition 1.4.16]{kr}).

Let $dV$ be a $C^{\infty}$-volume form on $X$. 
For a singular hemitian line bundle $(F,h_{F})$ on $X$,
let $A^{2}(M,F,h_{F},dV)$ denote the Hilbert space of 
$L^{2}$ holomorphic sections of $F$ with respect to $h_{F}$ 
and $dV$. 
We may also assume that for any pseudoeffective
singular hemitian line bundle
$(L,h_{L})$ and the any point $x\in X$, 
there exists an interpolation operator
\[
I_{x} : A^{2}(x,K_{X}\otimes A\otimes L,dV^{-1}h_{A}h_{L},\delta_{x})
\longrightarrow  A^{2}(X,K_{X}\otimes A\otimes L,dV^{-1}h_{A}h_{L},dV)
\]
such that the operator norm of $I_{x}$ is bounded from above 
by a positive constant independent of  $x\in X$ and $(L,h_{L})$, 
where $\delta_{x}$ denotes the Dirac measure at $x$. 
This is certainly possible, if we take $A$ to be 
sufficiently ample.

In fact let $x$ be a point on $X$ and let $U$ be a local 
coordinate neighbourhood of $x$ which is biholomorphic to $\Delta^{n}$. 
Then by the successive use of Theorem 2.4, we find an interpolation operator
\[
I_{x}^{U}: A^{2}(x,K_{X}\otimes A\otimes L,dV^{-1}h_{A}h_{L},\delta_{x})
\longrightarrow  A^{2}(U,K_{X}\otimes A\otimes L,dV^{-1}h_{A}h_{L},dV)
\]
such that the operator norm of $I_{x}^{U}$ is bounded from above 
by a positive constant independent of $(L,h_{L})$.
Now we note that the curvature of $h_{A}h_{L}$ is bounded from below by 
by the K\"{a}hler form $\Theta_{A}$.
Let $\rho$ be a $C^{\infty}$-function on $X$ such that 
$\mbox{Supp}\,\,\rho \subset\subset U$, $0\leq \rho\leq 1$ and 
$\rho \equiv 1$ on a neighbourhood of $x$. 
Let $\sigma_{x}$ be an element of $A^{2}(x,K_{X}\otimes A\otimes L,dV^{-1}h_{A}h_{L},\delta_{x})$.
Then  replacing  $(A,h_{A})$ to be some sufficiently high positive multiplie, 
by the usual $L^{2}$-estimate, we may assume that we can solve the 
$\bar{\partial}$-equation 
\[
\bar{\partial}u = \bar{\partial}(\rho\cdot I^{U}_{x}\sigma_{x} )
\]
with 
\[
u (x) = 0
\]
so that 
\[
\int_{X}h_{A}h_{L}\mid u\mid^{2}
\leq C\cdot (dV^{-1}h_{A}h_{L})(\sigma_{x},\sigma_{x})
\]
holds for a positive constant $C$ independent of $(L,h_{L})$ and 
$\sigma_{x}$.
Then 
\[
\rho\cdot I^{U}_{x}\sigma_{x} - u
\in H^{0}(X,{\cal O}_{X}(K_{X}+A+L)\otimes{\cal I}(h_{L}))
\]
is an extension of $\sigma_{x}$. 
Since $X$ is compact, moving $x$ and $U$, this implies the assertion.
\begin{lemma}
Let $h$ be a canonical AZD of $K_{X}$ constructed as in the proof 
of Theorem 2.2.
Then the inclusion :
\[
{\cal I}(h^{m}) \subseteq {\cal I}(h_{m})
\]
holds for every $m\geq 0$.  
\end{lemma}
{\bf Proof}. 
We prove this lemma by inducition on $m$. 
If $m = 0$, then both sides are ${\cal O}_{X}$. 

Suppose that the inclusion 
\[
{\cal I}(h^{m-1}) \subseteq {\cal I}(h_{m-1})
\]
has been settled for some $m\geq 1$. 
Then we have that by the property of $A$ as above
\[
{\cal O}_{X}(K_{X}+ (A + (m-1)K_{X}))\otimes {\cal I}(h^{m-1})
\]
is a globally generated subsheaf of 
\[
{\cal O}_{X}(A + mK_{X})\otimes {\cal I}(h_{m-1})
\]
Hence by the definition of ${\cal I}(h_{m})$ we see that 
\[
{\cal I}(h^{m-1}) \subseteq {\cal I}(h_{m})
\]
holds.  In particular 
\[
{\cal I}(h^{m}) \subseteq {\cal I}(h_{m})
\]
holds. 
By the induction on $m$, this completes the proof of Lemma 3.1.
{\bf Q.E.D.} \vspace{5mm} \\
By the choice of $A$ and Lemma 3.1, $h_{m}$ is well defined 
for every $m\geq 0$. 
Now we shall make the above lemma quantitative. 
\begin{lemma}
There exists a positive constant $C$ such that 
\[
h_{m} \leq C^{m}\cdot h_{A}\cdot h^{m} 
\]
holds for every $m\geq 0$. 
\end{lemma}
{\bf Proof}. 
First if   $m = 0$, both sides are $h_{0}$. 
Suppose that for some $m \geq 1$, 
\[
h_{m-1} \leq C_{(m-1)}\cdot h_{A}\cdot h^{m-1}
\]
holds for some positive constant $C_{(m-1)}$.
Let $dV$ be a $C^{\infty}$ volume form on $X$. 
Let $C(1)$ be a positive constant such that 
\[
h \geq C(1)\cdot (dV)^{-1}
\] 
holds on $X$. 
Let us denote the Bergman kernel of 
$A +mK_{X}$ with respect to a singular hermitian metric 
$H$ on $A+mK_{X}$ and the volume form $dV$ by 
$K(A+mK_{X},H,dV)$. 
In this notation $K_{m}$ is expressed as 
$K(A+ mK_{X},h_{m-1}\otimes (dV)^{-1},dV)$.

Then we have that 
\begin{eqnarray*}
K(A+ mK_{X},h_{m-1}\otimes (dV)^{-1},dV) 
& \geq & C_{(m-1)}^{-1}\cdot K(A+mK_{X},h_{A}\cdot h^{m-1}\cdot (dV)^{-1},dV) 
\\
& \geq & C_{(m-1)}^{-1}\cdot C(2)\cdot (h_{A}\cdot h^{m-1}\cdot (dV)^{-1})^{-1} \\
& \geq & C_{(m-1)}^{-1}\cdot C(1)\cdot C(2)
\cdot (h_{A}\cdot h^{m})^{-1}
\end{eqnarray*}
hold for some positive constant $C(2)$. 
The first inequality follow from
the formula : 
\[
K(A+ mK_{X},h_{m-1}\otimes (dV)^{-1},dV)(x)
\]
\[
\hspace{15mm} = \sup \{ \mid\sigma\mid^{2}(x) ; 
\sigma\in\Gamma (X,{\cal O}_{X}(A+mK_{X})) ;  
\int_{X}h_{m-1}\mid \sigma\mid^{2} =1\} \,\, (x\in X)
\]
and the similar formula for $K(A+mK_{X},h_{A}\cdot h^{m-1}\cdot (dV)^{-1},dV)$.

The 2-nd inequality follows from the $L^{2}$-extension theorem  (Theorem 2.4)  
applied to the extension from a point to $X$ as in Section 3. 
Hence we may assume that $C(2)$ is independent of $m$. 
Now we can take $C$ to be the constant 
$C(1)^{-1}\cdot C(2)^{-1}$.
This completes the proof of Lemma 3.2. 
{\bf Q.E.D.}
\vspace{5mm} \\
\begin{lemma}
There exists a positive constant $\tilde{C}$ such that 
for every $m\geq 1$,
\[
h_{A}\cdot K_{m}\leq \tilde{C}^{m}(dV)^{m}
\]
holds.
\end{lemma}
{\bf Proof}. 
Let $p\in X$ be an arbitrary point. 
Let $(U,z_{1},\ldots ,z_{n})$ be a local cooordinate around $x$ 
such that 
\begin{enumerate}
\item $z_{1}(p) = \cdots = z_{n}(p) = 0$,
\item $U$ is biholomorphic to the open unit polydisk in 
$\mbox{\bf C}^{n}$ with center $O\in \mbox{\bf C}^{n}$ 
by the coordinate,
\item $z_{1},\ldots ,z_{n}$ are holomorphic on a neighbourhood 
of the closure of $U$,
\item there exists a holomorphic frame $\mbox{\bf e}$ of $A$ on the closure of $U$.
\end{enumerate}
We set 
\[
\Omega := (\sqrt{-1})^{\frac{n(n-1)}{2}}dz_{1}\wedge\cdots\wedge dz_{n}\wedge d\bar{z}_{1}\wedge
\cdots\wedge d\bar{z}_{n}.
\]
For every $m\geq 0$, we set 
\[
C_{m} := \sup_{x\in U}\frac{h_{A}\cdot K_{m}}{\Omega^{m}}.
\]
We note that for any $x\in X$
\[
K_{m}(x) = \sup\{ \mid\phi\mid^{2}(x) ; 
\phi\in\Gamma (X,{\cal O}_{X}(A+mK_{X})), 
\int_{X}h_{m-1}\mid\phi\mid^{2} = 1\}
\]
holds.
Let $\phi_{0}$ be the element of $\Gamma (X,{\cal O}_{X}(A+mK_{X}))$ 
such that 
\[
K_{m}(x) = \mid\phi_{0}\mid^{2}(x)
\]
and 
\[
\int_{X}h_{m-1}\mid\phi_{0}\mid^{2} = 1.
\]
Then there exists a holomorphic function $f$ on $U$ such that 
\[
\phi_{0}\mid_{U} = f\cdot (dz_{1}\wedge\cdots dz_{n})^{m}\cdot\mbox{\bf e}
\]
holds.
Then 
\[
\int_{U}h_{A}\mid\phi_{0}\mid^{2}\Omega^{-(m-1)}
= \int_{U}\mid f\mid^{2}h_{A}(\mbox{\bf e},\mbox{\bf e})\Omega
\]
holds.  
On the other hand by the definition of $C_{m-1}$ we see that
\[
\int_{U}h_{A}\mid\phi_{0}\mid^{2}\Omega^{-(m-1)} 
\leq C_{m-1}\int_{U}h_{m-1}\mid\phi_{0}\mid^{2}
\leq C_{m-1}
\]
hold.
Combining above inequalities we have that
\[
\int_{U}\mid f\mid^{2}h_{A}(\mbox{\bf e},\mbox{\bf e})\,\Omega
\leq C_{m-1}
\]
holds.
Let $0 < \delta << 1$ be a sufficiently small number. 
Let $U_{\delta}$ be the inverse image of 
\[
\{ (y_{1},\ldots ,y_{n})\in \mbox{\bf C}^{n}\mid 
\mid y_{i}\mid < 1 - \delta \}
\]
by the coordinate $(z_{1},\ldots ,z_{n})$.

Then by the subharmonicity of $\mid f\mid^{2}$, there exists a positive constant $C_{\delta}$ independent of $m$ such that 
\[
\mid f(x)\mid^{2} \leq C_{\delta}\cdot C_{m-1}
\]
holds for every $x\in U_{\delta}$.
Then we have that 
\[
K_{m}(x) \leq C_{\delta}\cdot C_{m-1}\cdot\mid\mbox{\bf e}\mid^{2}\otimes\Omega^{m}(x)
\]
holds for every $x\in U_{\delta}$. 
Hence moving $p$, by the compactness of $X$ we see that 
there exists a positive constant $\tilde{C}$ such that 
\[
h_{A}\cdot K_{m} \leq \tilde{C}^{m}\cdot (dV)^{m}
\]
holds on $X$.
This completes the proof of Lemma 3.3. \vspace{10mm} Q.E.D. \\
\vspace{5mm} \\ 
By Lemma 3.2 and Lemma 3.3 
\[
K_{\infty} : = \mbox{the uppersemicontinuous envelope of} 
\,\,\,\, \limsup_{m\rightarrow\infty}
\sqrt[m]{K_{m}}
\]
is a well defined volume form on $X$ which does not vanish 
outside of a set of measure $0$. 
We set 
\[
h_{\infty} : = \frac{1}{K_{\infty}}. 
\]
Then by Lemma 3.2, we see that 
\[
h_{\infty} \leq C\cdot h
\]
holds. 
By the definition of $h_{\infty}$, it is clear that 
the curvature $\Theta_{h_{\infty}}$ is semipositive 
in the sense of current. 
Hence 
by the construction of $h$ (see the proof of Theorem 2.2), 
we see that there exists a positive constant $C^{\prime}$ 
such that the opposite estimate :  
\[
h_{\infty} \geq C^{\prime}\cdot h
\]
holds. 

Hence we have the following theorem.
\begin{theorem}
Let $h_{\infty}$ be the above singular hermitian metric on $K_{X}$.
Then $h_{\infty}$ is an AZD of $K_{X}$.  
\end{theorem}

\section{Proof of Theorem 1.1}
Now we shall prove Theorem 1.1. 
Let $\pi : X \longrightarrow \Delta$ be 
a smooth projective family of projective varieties 
as in Theorem 1.1.
We set $X_{t} := \pi^{-1}(t)$. 
If for every $t\in \Delta$, $K_{X_{t}}$ is not pseudoeffective,
then $P_{m}(X_{t}) = 0$ holds for every $t\in \Delta$ and 
every $m\geq 1$.
Hence in this case there is nothing to prove. 

Now we shall assume that for some $t_{0}\in \Delta$, say 
$t_{0} = 0$, $K_{X_{0}}$ is pseudoeffective. 
Shrinking $\Delta$, if necessary, we may 
assume that there exists an ample line bundle $A$ on $X$ such that 
for any pseudoeffective singular hermitian line bundle $(L,h_{L})$ 
\[
{\cal O}_{X}(A + L)\otimes{\cal I}(h_{L})
\]
and 
\[
{\cal O}_{X}(K_{X}+A+L)\otimes{\cal I}(h_{L})
\]
are globaly generated and for every $t\in \Delta$ and for any pseudoeffective singular hermitian 
line bundle $(L_{t},h_{L_{t}})$ on $X_{t}$,  
\[
{\cal O}_{X_{t}}(A\mid_{X_{t}} + L_{t})\otimes{\cal I}(h_{L_{t}})
\]
and 
\[
{\cal O}_{X_{t}}(K_{X_{t}}+A\mid_{X_{t}}+L_{t})\otimes{\cal I}(h_{L_{t}})
\]
are globaly generated.
Let $h_{A}$ be a $C^{\infty}$-hermitian metric on $A$ 
such that $\Theta_{h_{A}}$ is a K\"{a}hler form on $X$. 
We set 
\[
h_{0} := h_{A} 
\]
and 
\[
h_{0,t} := h_{A}\mid_{X_{t}} \,\,\,\, (t\in \Delta ).
\]
As in Section 3, inductively we shall define the sequences of 
singular hemitian metrics $\{ h_{m}\}$  on $X$ and
$\{ h_{m,t}\}$ on $X_{t} (t\in\Delta )$.
In this case $X$ is noncompact, but the construction 
works as in Section 3.   
However we should note that we do not know the psuedoeffectivity 
of $K_{X}$ or $K_{X_{t}}(t\in \Delta^{*})$ apriori. 
Hence at this stage $h_{m}$ and $h_{m,t}(t\in \Delta^{*})$
are {\bf not well defined} for $m\geq 2$. 

But by the $L^{2}$-extension theorem (Theorem 2.4 or 
\cite{m}) (as in the proof of Lemma 4.1),  we have the following lemma. 
\begin{lemma}
Let $t$ be a point on $\Delta$. 
Suppose that $h_{m-1}, h_{m-1,t}$ 
have been defined and 
\[
{\cal I}(h_{m-1,t})\subseteq {\cal I}(h_{m-1}\mid_{X_{t}})
\]
holds on $X_{t}$.
Then  every element of 
\[
H^{0}(X_{t},{\cal O}_{X_{t}}(K_{X_{t}}+ A\mid_{X_{t}}+(m-1)K_{X_{t}})
\otimes {\cal I}(h_{m-1,t}))
\]
extends to an element of 
\[
H^{0}(X,{\cal O}_{X}(K_{X}+ A+(m-1)K_{X})\otimes {\cal I}(h_{m-1})).
\]

\end{lemma}
{\bf Proof}.
We note that since Theorem 2.4 is stated for Stein manifold. 
Hence  we cannot apply Theorem 2.4 directly.  
Let $U$ be a Zariski open Stein subset of $X$ such that $X_{t}\cap U$ 
is nonempty and $K_{X}\mid_{U}$ is holomorphically trivial. 
Then for every element $\sigma_{t}$ of 
\[
H^{0}(X_{t},{\cal O}_{X_{t}}(K_{X_{t}}+ A\mid_{X_{t}}+(m-1)K_{X_{t}})
\otimes {\cal I}(h_{m-1,t})), 
\]
$\sigma_{t}\mid_{U\cap X_{t}}$ extends to an element $\sigma$ of 
\[
H^{0}(U,{\cal O}_{X}(K_{X}+ A+(m-1)K_{X})\otimes {\cal I}(h_{m-1})).
\]
with the $L^{2}$-condition 
\[
\int_{U}h_{m-1}\mid\sigma\mid^{2} < \infty .
\]
But this $L^{2}$-condition implies that $\sigma$ extends to a section of 
\[
H^{0}(X,{\cal O}_{X}(K_{X}+ A+(m-1)K_{X})\otimes {\cal I}(h_{m-1})).
\]
\vspace{5mm} {\bf Q.E.D.}\\

Since $K_{X_{0}}$ is pseudoeffective, using  Lemma 4.1, 
we have that $h_{m}$ is well defined for every $m\geq 0$ and 
 the inclusion ${\cal I}(h_{m,0})\subset {\cal I}(h_{m}\mid_{X_{0}})$
holds for every $m\geq 0$ inductively. 
Hence for every $m\geq 1$ we have a proper Zariski closed subset 
$S_{m}$ in $\Delta$ such that 
for every $t\in \Delta -S_{m}$,  $h_{m,t}$ is well defined.
In particular $K_{X_{t}}$ is pseudoeffective 
for every $t\in \Delta - \cup_{m\geq 1}S_{m}$. 
This implies that $K_{X_{t}}$ is pseudoeffective for every $t\in\Delta$.
Then using Lemma 4.1 and Lemma 3.1, by  induction on $m$, we have the following lemma.  
\begin{lemma}
$ h_{m}, h_{m,t} (t\in \Delta )$ are well defined 
for every $m\geq 0$ and $t\in \Delta$.
 And 
\[
{\cal I}(h_{m,t})\subseteq {\cal I}(h_{m}\mid_{X_{t}}) 
\]
holds for every $m\geq 0$ and $t\in \Delta$. 
\end{lemma}
We define 
\[
h_{\infty} := \mbox{the lowersemicontinuous envelope of} \,\, 
\liminf_{m\rightarrow \infty}\sqrt[m]{h_{m}}
\]
and 
\[
h_{\infty ,t} := \mbox{the lowersemicontinuous envelope of} \,\, 
\liminf_{m\rightarrow \infty}\sqrt[m]{h_{m,t}}.
\]
Although $X$ is noncompact, the argument in Section 3 
is still valid. 
In fact since $X$ admits a continuous plurisubharmonic 
exhaustion function, Lemma 3.2 holds in this case, if 
we restrict the family to a relatively compact subset
 of $\Delta$. 
Also  Lemma 3.3 is valid on every relatively 
compact subset of $X$, since the proof is local.   
Hence $h_{\infty}$, $h_{\infty ,t}$ are well defined AZD's
on $K_{X}$ and $K_{X_{t}}$ respectively. 

Again by the $L^{2}$-extension theorem (Theorem 2.4), as  Lemma 3.2, we have the following lemma.  
\begin{lemma}
For every $t$, there exists a positive constant $C$ such that 
\[
h_{m}\mid_{X_{t}} \leq C^{m}\cdot h_{m,t} 
\]
holds for every $m\geq 0$.
In particular 
\[
h_{\infty}\mid_{X_{t}} \leq C\cdot h_{\infty ,t}
\]
holds. 
\end{lemma}
Here we have used the fact that 
for every 
\[
\sigma_{t}\in A^{2}(X_{t},K_{X_{t}}\otimes K_{X}^{\otimes (m-1)},h_{m-1}\mid_{X_{t}}),
\]
there exists an element of 
\[ 
\sigma \in A^{2}(X_{t},K_{X}^{\otimes m},h_{m-1})
\]
such that $\sigma\mid_{X_{t}} = \sigma_{t}$ and  
\[
\parallel \sigma\parallel \leq C_{t}\parallel \sigma\parallel ,
\]
where $C_{t}$ is a positive constant depending only on $t$ (if $t= 0$, we may 
take $C_{0}$ to be $2\sqrt{1620\pi }$ by Theorem 2.4). 
Here we note that the $L^{2}$-spaces above are determined without specifying 
volume forms. 

Now by Lemma 4.3 and Theorem 3.1, we see that 
$h_{\infty}\mid_{X_{t}}$ is an AZD of $K_{X_{t}}$.
Then by the $L^{2}$-extension theorem (Theorem 2.4, 
see also the proof of Lemma 4.1), we see that for every $m\geq 1$, 
every element of 
\[
H^{0}(X_{t},{\cal O}_{X_{t}}(mK_{X_{t}})\otimes{\cal I}(h_{\infty }^{m-1}\mid_{X_{t}}))
\simeq H^{0}(X_{t},{\cal O}_{X_{t}}(mK_{X_{t}}))
\]
extends to an element of 
\[
H^{0}(X,{\cal O}_{X}(K_{X}+(m-1)K_{X})\otimes
{\cal I}(h_{\infty}^{m-1})).
\]
Hence we see that $P_{m}(X_{t})$ is lowersemicontinuous.
By the upper semicontinuity theorem for cohomologies,
we see that $P_{m}(X_{t})$ is independent of $t\in \Delta$. 
This completes the proof of Theorem 1.1. 

Theorem 1.2 has already been proved in the course of the above proof.

Author's address\\
Hajime Tsuji\\
Department of Mathematics\\
Tokyo Institute of Technology\\
2-12-1 Ohokayama, Megro 152-8551\\
Japan \\
e-mail address: tsuji@math.titech.ac.jp
\end{document}